\documentclass[16pt,reqno,oneside]{etds}
\usepackage{amsfonts,amssymb,latexsym,xspace,epsfig,graphics,color}
\usepackage{amsmath,enumerate,stmaryrd,xy}
\usepackage{bbm}
\usepackage[numbers]{natbib}
\usepackage{subfigure}

\newcommand{\E}{\mathbb{E}}

\newcommand{\N}{\ensuremath{\mathbb{N}}}

\newcommand{\R}{\ensuremath{\mathbb{R}}}     
\newcommand{\Z}{\ensuremath{\mathbb{Z}}}    
\renewcommand{\P}{\ensuremath{\mathbb{P}}}

\def\eqd{\,{\buildrel d \over =}\,}

\usepackage{algorithmic}
\usepackage{algorithm, caption}

\newtheorem{lemma}{Lemma}

\newtheorem{obs}{Observation}

\newtheorem{theo}{Theorem}

\newtheorem{coro}{Corollary}

\DeclareGraphicsExtensions{.pdf,.gif,.jpg}

\begin{document} 
\ETDS{0}{0}{0}{0} 
\runningheads{S. Gallo and D. Y. Takahashi}{Attractive Regular Stochastic Chains} 
\title{Attractive regular stochastic chains:\\ perfect simulation and phase transition\footnotetext{This work is part of USP project ``Mathematics, computation, language and the brain'' and partially elaborated during the meeting ``Jorma's Razor 2" supported by CAPES grant AUXPE-PAE-598/2011. SG was supported by FAPESP grant 2009/09809-1. DYT was partially supported by FAPESP grant 2008/08171-0 and Pew Latin American Fellowship.}} 
\author{SANDRO GALLO\affil{1}\, and DANIEL Y. TAKAHASHI\affil{2}\, } 
\address{\affilnum{1}\ Instituto de Matem\'atica, Universidade Federal de Rio de Janeiro, Brasil.\\ 
\email{sandro@im.ufrj.br}\\
\affilnum{2}\ Neuroscience Institute and Psychology Department, Princeton University, USA.\\ 
\email{takahashiyd@gmail.com}} \\
\recd{...}



%
%
%
%
%
%

\begin{abstract}
We prove that uniqueness of the stationary chain, or equivalently, of the $g$-measure, compatible with an attractive regular probability kernel is equivalent to either one of the following two assertions for this chain: (1) it is a finitary coding of an i.i.d. process with countable alphabet, (2) the concentration of measure holds at exponential rate. We show in particular that if a stationary chain is uniquely defined by a kernel that is continuous and attractive, then this chain can be sampled using a coupling-from-the-past algorithm. For the original Bramson-Kalikow model we further prove that there exists a unique compatible chain if and only if the chain is a finitary coding of a finite alphabet i.i.d. process. Finally, we obtain some partial results on conditions for phase transition for general chains of infinite order.
\end{abstract}

\section{Introduction}

In this work we consider chains of infinite order, or equivalently $g$-measures, on finite alphabet. These are processes specified by kernels of transition probabilities that can depend on the whole past and include as special cases the finite order Markov chains, stochastic models that exhibit phase transitions \citep{hulse/2006,berger/hoffman/sidoravivius/2005, bramson/kalikow/1993}, and models that are not Gibbsian \cite{fernandez/gallo/maillard/2011}. 
An important question for this class of models is ``what properties distinguish kernels exhibiting phase transition from kernels satisfying uniqueness?''
This work gives necessary and sufficient conditions for the existence of phase transition for an important class of chains of infinite order, namely for the attractive regular chains.

A probability kernel is called regular when it is strongly non-null and continuous with respect to the past. The regularity of the kernel guarantees the existence of at least one chain compatible with the kernel. Attractiveness means that the transition probabilities exhibit a monotonicity property and it is analogous to the attractiveness of the specifications considered in statistical mechanics \cite{preston/1976}.  

It was quite unexpected when \citet{bramson/kalikow/1993} showed an example of a family of regular  and attractive kernels with more than one compatible chain. Some interesting results are based on this example (\emph{BK example} in the sequel). For instance, \citet{quas/1996}  used the BK example to construct a $C^1$ expanding map of the circle which preserves Lebesgue measure and such that the system is ergodic, but not weak-mixing. Also using the BK example, \citet{stenflo/2001} showed a counterexample to a conjecture raised by \citet{karlin/1953}. 
\citet{lacroix/2000}  obtained some simplification to the \citet{bramson/kalikow/1993} proof of phase transition and  \citet{hulse/2006} showed a different example of a regular and attractive kernel exhibiting phase transition using similar ideas of  proofs as in \cite{lacroix/2000}.  To the best of our knowledge, \citet{berger/hoffman/sidoravivius/2005} exhibited the only non-attractive example of phase transition that is not based on the BK example. Despite the importance of these works, none of them give general sufficient conditions for the existence of phase transition, even for the special (but important) case of attractive kernels.

 In the present work, we prove that for attractive regular kernels, uniqueness of the stationary chain compatible with a kernel is equivalent to either one of the following two assertions: (1) the compatible chain is a finitary coding of a countable alphabet i.i.d. process, (2) the concentration of measure holds at exponential rate. Condition (1) means that the compatible chain is a factor of a countable alphabet i.i.d. process whose mapping depends almost surely on a finite number of coordinates with respect to the measure of the i.i.d. process. Condition (2) means that phase transition yields loss of  ``good'' concentration of measure. 
We show that the regularity of the kernel is essential for the three equivalences to hold and, in general, cannot be relaxed. We also obtain some partial results for non-regular and non-attractive cases, which are of independent interest.

The main ingredient for the proof of the existence of a finitary coding is the proof that uniqueness of the compatible chain for regular and attractive kernels is equivalent to the existence of a coupling-from-the-past (CFTP) perfect simulation algorithm. This class of simulation algorithms was first introduced  for Markov chains by \citet{propp/wilson/1996} and then generalized to several other stochastic models.  Despite its simplicity, our algorithm can generate samples of continuous and attractive chains under regimes in which it was previously not known to be possible.   

Finally, from a coding point of view, it is interesting to have a finitary coding from a finite alphabet i.i.d. process. We show that this is possible in the original  \citet{bramson/kalikow/1993} example if and only if there is a unique compatible chain, that is, choosing the parameters of the model in such a way that uniqueness holds.

It is worth mentioning that \citet{steif/vanderberg/1999} obtained similar results for a certain class of $d$-dimensional Probabilistic Cellular Automata.  Specifically, they proved that if a random field is obtained as invariant measures of monotonic, exponentially ergodic probabilistic cellular automata,  then the random field is a finitary coding of a finite alphabet i.i.d. random field. They used this result to show that there exists a finitary coding from a finite alphabet  i.i.d. process to the plus phase of a ferromagnetic Ising model strictly below the critical temperature. They also proved that for the plus phase of an Ising model strictly above the critical temperature, there is no finitary coding from a finite alphabet i.i.d. process. The situation for the critical temperature remains open. 
For the BK  example, it is not clear what is the quantity equivalent to the statistical mechanics notion of temperature, and therefore we cannot make a one-to-one comparison with Theorem 1.1 in \cite{steif/vanderberg/1999}. Nevertheless, we observe that our result for the original BK example does not have a priori restriction on the parameters of the model or on the loss of memory of the process. 

This article is organized as follows. In Section \ref{sec:def} we introduce the  notation, definitions, the necessary background, and two important examples. In Section \ref{sec:results} we state the main results. In Section \ref{sec:attractivesampler} we introduce the Attractive Sampler, which is used to prove the theorems of Section \ref{sec:results}. Finally we prove the results in Section \ref{sec:proof}.

\section{Notation, standard definitions, and examples}\label{sec:def}

\paragraph{{\bf Notation}}
For any set $\mathcal{U}$ we denote the sets of bi-infinite, right-infinite and finite sequences of symbols of $\mathcal{U}$ by $\,\mathcal{U}^{-\mathbb{Z}}=\mathcal{U}^{\{\ldots,2,1,0,-1,-2,\ldots\}}$, $\mathcal{U}^{-\N}=\mathcal{U}^{\{-1,-2,\ldots\}}$ and $\mathcal{U}^\star = \bigcup_{j \geq1} \mathcal{U}^{\{-1, \ldots,-j\}}$, respectively. The elements of these sets will be denoted, respectively, ${\bf u}=\ldots u_{2}u_{1}u_{0}u_{-1}u_{-2}\ldots$, $\underline{u}=u_{-1}u_{-2}\ldots$ and $u_{-k}^{-1}=u_{-1}u_{-2}\ldots u_{-k}$ for any $1\leq k<+\infty$. We also use the notation $\,\mathcal{U}^{-n}=\mathcal{U}^{\{-1,\ldots,-n\}}$ and $\,\mathcal{U}^{n} =\mathcal{U}^{\{n,\ldots,1\}} $for any $n\ge1$. Finally, if $u\in \mathcal{U}^{n}$ and $v\in \mathcal{U}^{m}$ we will use the notation $uv$ for the string of length $n+m$ belonging to $\mathcal{U}^{n+m}$, obtained by concatenating both strings. This notation is extended to the case where one of the string is a right-infinite sequence, for example $v\underline{u}$.
Note that we are using the convention that the past (negative indices) is on the right and the future (positive indices) is on the left.

In the present article $\mathcal{U}$ is some Polish space and $A$ is the finite ordered set $\{1,2,\ldots,s\}$ unless specified. $A$ is called \emph{alphabet.} We define a partial order on $A^{-\mathbb{N}}$ by saying that $\underline{a}\leq \underline{b}$ whenever $a_{-i}\leq b_{-i}$ for every $i\ge1$. In $A^{-\mathbb{N}}$, the maximal element is $\underline{s}:=sss\ldots $ and the minimal element is $\underline{1}:=111\ldots$.

\paragraph{{\bf Chains of infinite order}}A \emph{probability kernel}, or simply a \emph{kernel},  $P$ on the alphabet $A$ is a function  
\begin{equation*}
\begin{array}{cccc}
P:&A\times A^{-\mathbb{N}}&\rightarrow& [0,1]\\
&(a,\underline{x})&\mapsto&P(a|\underline{x})
\end{array}
\end{equation*}
such that
\[
\sum_{a\in A}P(a|\underline{x})=1\,\,,\,\,\,\,\,\,\forall \underline{x}\in A^{-\mathbb{N}}.
\] 
We say that a stationary stochastic chain ${\bf X} = \{X_{j}\}_{j\in\Z}$ (of stationary law $\mu$) on $A$ is \emph{compatible} with a kernel $P$ if the later is a regular version of the conditional probabilities of the former, that is,
\begin{equation*}\label{compa}
\mu(X_{0}=a|X_{-\infty}^{-1}=\underline{x})=P(a|\underline{x})
\end{equation*}
for every $a\in A$ and $\mu$-a.e. $\underline{x}$ in $A^{-\mathbb{N}}$. When there is more than one stationary chain compatible with $P$, we say that there is \emph{phase transition}, otherwise we say that the chain is \emph{unique}. We follow the Harris nomenclature \cite{harris/1955} and call \emph{chains of infinite order} the chains compatible with kernels. 
They were first introduced in \citep{onicescu/mihoc/1935} under the name \emph{cha\^ines \`a liaisons compl\`etes}. The existence of an invariant measure for these chains was first studied by \citet{doeblin/fortet/1937}. In ergodic theory, it was rediscovered by \citet{keane/1972} who gave the name $g$-function to the kernel and called $g$-measure the chain compatible with the kernel. For a comprehensive historical account and recent developments we refer the reader to \cite{fernandez/maillard/2005}.

\paragraph{{\bf Non-nullness, continuity rate, oscillations and attractiveness}}We say that a kernel $P$ is  \emph{strongly non-null} if 
\begin{equation*}
\inf_{a \in A, \underline{x}\in A^{-\mathbb{N}}}P(a|\underline{x}) > 0.
\end{equation*}
The \emph{continuity rate} (or \emph{variation}) \emph{of order} $k$ of a kernel $P$ is 
\begin{equation*}
\textrm{var}_{k}:=\sup_{b\in A}\sup_{a_{-k}^{-1}\in A^{-k}}\sup_{\underline{x},\underline{y}\in A^{-\mathbb{N}}}|P(b|a_{-k}^{-1}\underline{x})-P(b|a_{-k}^{-1}\underline{y})|. 
\end{equation*}
We say that $P$ is continuous if $\lim_{k\rightarrow \infty}\textrm{var}_{k}=0$. A compactness argument shows that if the kernel is continuous, at least one compatible stationary chain exists (see for example \cite{keane/1972}). If $P$ is strongly non-null and continuous, we say that $P$ is a \emph{regular kernel}. 

Another characterization of kernels is given by the \emph{oscillation rate}:
\[
\textrm{osc}_{n}:=\sum_{a\in A}\textrm{osc}_{n}(a)
\]
where
\[\textrm{osc}_{n}(a):=\sup\{|P(a|\underline{x})-P(a|\underline{y})|:\underline{x},\underline{y}\in A^{-\mathbb{N}}\,\,x_{-i}=y_{-i}\,\forall\,i\neq n\}.
\] 
The sequences $\{\textrm{var}_{k}\}_{k\ge0}$ and $\{\textrm{osc}_{k}\}_{k\ge0}$ are related to the uniqueness of the compatible stationary chain as we will see in the examples below.

Finally, we say that a kernel $P$ on $A$ is \emph{attractive} if for all $a\in A$ the value of $\sum_{j \geq a}P(j|\underline{x})$ is increasing on $\underline{x} \in A^{-\N}$. \\

Let us give two important examples taken from the literature, which we will repeatedly use in the sequel to illustrate our assertions.

\paragraph{{\bf Binary auto-regressive models}}

These models are extensively used in the statistical literature \cite{mccullagh/nelder/1983}, and are defined through the following parameters: a continuously differentiable and increasing function $\psi:\mathbb{R}\rightarrow]0,1[$ such that $\psi(t)+\psi(-t) = 1$,  a summable sequence of non-negative real numbers $\{\xi_{n}\}_{n\geq1}$, and a non-negative real parameter $\gamma\geq0$. Consider the class of kernels $P$ on the alphabet $\{-1,+1\}$ such that
\[
P(a|\underline{x}):=\psi\left(a\sum_{n\geq1}\xi_{n}x_{-n}+a\gamma\right).
\]
 A straightforward computation shows that these kernels are attractive and regular.
Moreover, if $\psi$ is Lipschitz continuous, then, one immediately obtains $\textrm{var}_{k}\leq C\sum_{n>k}\xi_{n}$ and $\textrm{osc}_{n}\leq C\xi_{n}$ for some positive constant $C$. In this case, possible criteria for uniqueness are $\xi_{n}=c n^{-\alpha}$ with $\alpha>3/2$ and $c$ any constant (\citet{johansson/oberg/2003}), or $\textrm{osc}_{n} \leq C\sum_{n\geq1}\xi_{n}<1$ (\citet{fernandez/maillard/2005}).

An important example of  binary auto-regressive models is when $\psi(t)=e^{-t}(e^{-t}+e^{t})$. The resulting kernel is called \emph{logit model} in the statistics literature, and \emph{one-sided $1$-dimensional long-range Ising model} in statistical physics literature. For instance, \citet{hulse/2006} used this model to give an example of phase transition in chains of infinite order.

\paragraph{{\bf The example of \citet{bramson/kalikow/1993}}}

Consider an increasing function $\phi:[-1,+1]\rightarrow]0,1[$ such that $\phi(t)+\phi(-t) = 1$, an increasing sequence of odd positive integers $\{m_{j}\}_{j\geq1}$ and a sequence $\{\lambda_{j}\}_{j\geq1}$ such that $\lambda_{j}\geq0$ and $\sum_{j\geq1}\lambda_{j}=1$. We call the \emph{BK example}  \citep{friedli/2010} the class of kernels defined on alphabet $A=\{-1,+1\}$ by
\begin{equation*}\label{eq:kernelBK}
P(+1|\underline{x})=\sum_{j\geq1}\lambda_{j}\phi\left(\frac{1}{m_{j}}\sum_{i=1}^{m_{j}}x_{-i}\right).
\end{equation*}
The kernels of this class are attractive and regular.  Attractiveness and strong non-nullness follow directly from the definition of $\phi$ and simple calculations yield $\textrm{var}_{n}\leq\sum_{\{j:m_{j}>n\}}\lambda_{j}$,
showing that $P$ is continuous. When $\phi(t)=(1-\epsilon){\bf 1}\{t>0\}+\epsilon{\bf 1}\{t<0\}$ for some $\epsilon\in(0,1/2)$ and $\lambda_k = cr^{k}$ for some $r \in (2/3,1)$ and $c=(1-r)/r$, we call this model the \emph{original BK example}, as it is precisely the model introduced in  \cite{bramson/kalikow/1993}, where it is proved that, taking the sequence $\{m_{k}\}_{k\ge1}$ increasing sufficiently fast, the kernel $P$ exhibits phase transition.

\paragraph{{\bf Maximum and minimum phases for attractive kernels}}

Define, for any $\underline{x}\in A^{-\mathbb{N}}$, the fixed past chain ${\bf X}^{\underline{x}}$ by
\begin{equation*} \label{eq:nonstationary}
X^{\underline{x}}_{n}=\left\{
\begin{array}{lll}
x_{n}&\textrm{if}&n\leq 0\\
a&\textrm{with probability}&P(a|X^{\underline{x}}_{n-1}\ldots X^{\underline{x}}_{1}\underline{x})\,\,\,\textrm{otherwise}.
\end{array}
\right.
\end{equation*}
${\bf X}^{\underline{x}}$ is the non-stationary chain obtained by fixing the past $\underline{x}$ from time $0$ to $-\infty$, and ``running $P$ from this past''. 
For any $\underline{x}$, $j\geq1$, and $n\in\mathbb{Z}$, let ${\bf X}^{\underline{x},j}$ be the process defined by $X_{n}^{\underline{x},j}:=X_{n+j}^{\underline{x}}$. The attractiveness of $P$ implies that for $\underline{x}=\underline{1}$ and $\underline{s}$, the sequence of processes $\{{\bf X}^{\underline{x},j}\}_{j\geq1}$ is stochastically non-decreasing and non-increasing, respectively \cite{hulse/1991}, and therefore, the weak limits 
\begin{equation*}\label{attract}
{\bf X}^{\min}:=\lim_{j \rightarrow+\infty}{\bf X}^{\underline{1},j}\,\,\textrm{and }\,\,{\bf X}^{\max}:=\lim_{j\rightarrow+\infty}{\bf X}^{\underline{s},j}
\end{equation*}
exist and are stationary. If $P$ is continuous, then ${\bf X}^{\min}$ and ${\bf X}^{\max}$ are compatible with $P$. We call ${\bf X}^{\min}$ the \emph{minimum phase} and ${\bf X}^{\max}$ the \emph{maximum phase}. 

\paragraph{{\bf Finitary process and B-process}} Let  $T_{\mathcal{U}}$ and $T_{A}$ be the shift operators that act respectively on $\mathcal{U}^{\mathbb{Z}}$ and $A^{\mathbb{Z}}$ shifting the sequences by one unit: $T_{\mathcal{U}}({\bf u})=\{u_{i+1}\}_{i\in\mathbb{Z}}$ and $T_{A}({\bf a})=\{a_{i+1}\}_{i\in\mathbb{Z}}$. 
A stationary process ${\bf X}$ (with stationary law $\mu$) on the alphabet $A$ is a \emph{stationary coding} of a stationary process ${\bf U}$ (with stationary law $\mathbb{P}$) on $\mathcal{U}$ if there exists a  measurable function $\Phi:\mathcal{U}^{\mathbb{Z}}\rightarrow A^{\mathbb{Z}}$ which is \emph{translation equivariant} (that is $\Phi(T_{\mathcal{U}}({\bf U}))=T_{A}\Phi({\bf U})$) and such that $\mu=\mathbb{P}\circ \Phi^{-1}$.  We follow the nomenclature given in \cite{shields/1996} and call \emph{B-process} a process that is a stationary coding  of an i.i.d. process. A stationary coding is called \emph{finitary coding} if there exists \emph{stopping times} $\theta_1:\mathcal{U}^{\Z}\rightarrow \mathbb{N}\cup \{\infty\}$  and $\theta_2:\mathcal{U}^{\Z}\rightarrow \mathbb{N}\cup \{\infty\}$, both $\mathbb{P}$-a.s. finite, such that 
\begin{equation}\label{eq:checkable}
[\Phi({\bf U})]_{0}=[\Phi({\bf V})]_{0}\,\,\,\textrm{whenever}\,\,\,U_{-\theta_1({\bf U})}^{+\theta_2({\bf U})}=V_{-\theta_1({\bf U})}^{+\theta_2({\bf U})}.
\end{equation}
This last assumption means that the event $\{\theta_1({\bf U})=k\} \cap \{\theta_2({\bf U})=l\}$ is $\mathcal{F}(U_{-k}^{+l})$-measurable. In other words, the stopping times $\theta_1({\bf U}), \theta_2({\bf U})$ are checkable looking only at an a.s. finite number of $U_{i}$'s. We call \emph{finitary processes} (FP) the processes that are finitary coding of an i.i.d. process. 
We will often use the simplified notation $\theta_1$ and $\theta_2$ respectively for $\theta_1({\bf U})$ and $\theta_2({\bf U})$.

The notion of stationary coding comes from ergodic theory, and has a one-side analogue in the literature of stochastic processes, called the \emph{coupling-from-the past} algorithm (\emph{CFTP algorithm} in the sequel). Such algorithms, which were first introduced in \cite{propp/wilson/1996} for Markov chains, aim to construct the function $\Phi$ using the kernel $P$ and an a.s. finite number of past values of an i.i.d. process ${\bf U}$. If a CFTP algorithm is \emph{feasible} for a given kernel $P$, then the constructed stationary measure $\mu$ is a FP, because it is a particular finitary coding of ${\bf U}$, for which we can take $\theta_2=0$ and $\theta_1$ finite $\mathbb{P}$-a.s. For simplicity, in this case, we use the notation $\theta:=\theta_{1}$ and say that $\theta$ is the stopping time of CFTP algorithm.

In the literature, sometimes a process is called finitary only if the set $\mathcal{U}$ is finite. We do not assume this. In the special case of  $\mathcal{U}$ being finite (or countable), we say that the process is a \emph{finitary coding of a finite (countable) alphabet i.i.d. process}.

\paragraph{{\bf Exponential concentration of measure}}

Let $f: A^n \to \R$ be measurable. Define $\delta_j f = \sup \{|f(a^{n}_1) - f(b^n_{1})|: a_i = b_i,  \forall i\neq j\}$ and let $\delta f$ be the vector with $j$-th coordinate given by $\delta_j f$.  We say that the \emph{concentration of measure holds at exponential rate} for a stationary process ${\bf X}$  if, for all  $n > 0$, $\epsilon > 0$, and functions $f$, we have
\begin{equation}\label{eq:expoconcentration}
\P\left(\left|f(X^{n}_{1}) - \E[f(X^n_1)]\right| > \epsilon \right) \leq C\exp\left\{{-\frac{ g(\epsilon, \|\delta f\|_{\ell_1(\N)})}{\|\delta f \|^2_{\ell_2(\N)}}} \right\}
\end{equation}
where $g(\epsilon, \|\delta f\|_{\ell_1(\N)}) > 0$ and $C$ is a numerical constant.\\

In particular, the above inequality implies that for all $k \in \N$, $n \geq k$, and $h: A^k \to [0,1]$ we have
\begin{equation} \label{eq:exfreq}
\P\left(\left|\frac{1}{n-k+1}\sum_{j=0}^{n-k}h(X^{j+k}_{j+1}) - \E[h(X^k_1)]\right| > \epsilon\right) \leq Ce^{-ng_k(\epsilon)}
\end{equation}
where $g_k(\epsilon) > 0$ and $C$ is a numerical constant. We say that the \emph{ergodic theorem holds at exponential rate}  for a stationary process ${\bf X}$ if it satisfies (\ref{eq:exfreq}).

\section{Main results} \label{sec:results}

\vspace{1cm}
\begin{theo} \label{theo:attractive}
Let $P$ be an attractive regular kernel. The following are equivalent:
\begin{enumerate}
\item There exists a unique stationary chain compatible with $P$.
\item ${\bf X}^{\max}$ is a finitary coding of a countable alphabet i.i.d. process.
\item The concentration of measure holds at exponential rate for ${\bf X}^{\max}$.
\end{enumerate}
\end{theo}

From \cite{hulse/1991} we know that the maximum phase (\emph{resp.} the minimum phase) is always a B-process  regardless of being equal or different to the minimum phase (\emph{resp.} the maximum phase). Therefore, the fact to be a B-process does not  distinguish the presence or not of phase transition. Theorem \ref{theo:attractive} shows that, for a regular attractive kernel $P$, to be a finitary coding of an i.i.d process or to have concentration of measure at exponential rate distinguishe between existence or not of phase transition.

We now show that Theorem \ref{theo:attractive} is optimal in the class of attractive chains, in the sense that if we relax either continuity or strong non-nullness, we can find examples of stationary chains that are FP and with the ergodic theorem holding at exponential rate, although they are not uniquely determined by their conditional probabilities.

\paragraph{{\bf Relaxing the strong non-nullness assumption}}

The following example shows that in general we cannot relax the strong non-nullness condition. Before giving our example, we need some more definitions. Let $A=\{-1,+1\}$. For any $\underline{x}\in A^{-\mathbb{N}}$, we define
\[
\ell(\underline{x}):=\min\{j\geq0:x_{-j-1}=-1\}.
\]
When we look into the past in $\underline{x}$, $\ell(\underline{x})$ counts the number of $+1$  before finding the first $-1$. Let $+\underline{1}$ be the pasts such that $+\underline{1}_j = +1$ for all $j$. We use the convention that $\ell(+\underline{1})=\infty$.
Let $\{p_{i}\}_{i\geq0}$ be a monotonically decreasing sequence of $[0,1]$-valued real numbers, and let $p_{\infty}=\lim_{i\rightarrow+\infty}p_{i}$. The kernel $P$ is defined on $\{-1,+1\}$ by $P(-1|\underline{x})=p_{\ell(\underline{x})}$ for any $\underline{x}\neq +\underline{1}$ and $P(-1|+\underline{1}):=p_{\infty}$. It is clear that this example is attractive. It is also continuous. To see this, observe that
\begin{equation*}
\sup_{\underline{y},\underline{z}\in A^{-\mathbb{N}}}|P(a|a_{-k}^{-1}\underline{y})-P(a|a_{-k}^{-1}\underline{z})|=0
\end{equation*}
for any $a_{-k}^{-1}\in A^{k}$, except for $(+1)^{k}$. Hence,
\[
\sup_{a_{-k}^{-1}\in A^{k}}\sup_{\underline{y},\underline{z}\in A^{-\mathbb{N}}}|P(a|a_{-k}^{-1}\underline{y})-P(a|a_{-k}^{-1}\underline{z})|=\sup_{\underline{y},\underline{z}\in A^{-\mathbb{N}}}|P(a|(+1)^{k}\underline{y})-P(a|(+1)^{k}\underline{z})|,
\]
and thus we obtain that 
\[
\sup_{a_{-k}^{-1}\in A^{k}}\sup_{\underline{y},\underline{z}\in A^{-\mathbb{N}}}|P(a|a_{-k}^{-1}\underline{y})-P(a|a_{-k}^{-1}\underline{z})|=p_{k}-p_{\infty},
\]
which goes to $0$ by the definition of the sequence $\{p_{k}\}_{k\geq0}$. If $p_{\infty}=0$, the chain is not strongly non-null, and the degenerated chain with all symbols equal to  $+1$ is trivially stationary and compatible with $P$. If we further assume $\sum_{k\geq1}\prod_{i=0}^{k-1}(1-p_{i})<+\infty$, there exists another class of stationary chains compatible with the kernel $P$. It is the so-called \emph{renewal chains}, obtained by concatenating i.i.d. blocks of the form $(-1,+1,\ldots,+1,+1)$  of random length. These blocks have length $k+1$ with probability $\prod_{i=0}^{k-1}(1-p_{i})p_{k}$, and therefore, have finite expected length. The existence of several compatible chains is due to the fact that this kernel is \emph{not irreducible} \cite{cenac/chauvin/paccaut/pouyanne/2010}. Therefore, this chain has a degenerate type of phase transition. Nevertheless, the compatible chain with probability one of having $+1$ is obviously a finitary coding of an i.i.d. process with a concentration of measure at exponential rate. This shows that if the strong non-nullness assumption is removed, the existence of a finitary coding for the maximum phase does not imply uniqueness of the compatible chain.

\paragraph{{\bf Relaxing the continuity assumption}} For a discontinuous attractive kernel $P$, the maximum and minimum phases are always distinct and not consistent with the kernel $P$ \citep{hulse/1991}. Hence, strictly speaking, we don't have a phase transition where there is more than one chain compatible with $P$. In fact, this means that considering discontinuous attractive chains does not make much sense from the point of view of non-uniqueness. Nevertheless, we can still ask if the maximum and minimum phases of a discontinuous attractive kernel can be finitary codings of  i.i.d process. The example below shows that, in general, this could happen.

Let $\underline{a}\in  \{-1,1\}^{-\N}$ and let $-\underline{1}$ be the pasts such that $-\underline{1}_j = -1$ for all $j$. More generally, for any $\underline{z}\in\{-1,1\}^{-\N}$ (\emph{resp.} $u\in \{-1,1\}^{\star}$) we denote by $-\underline{z}$ (\emph{resp.} $-u$) the sequence of symbols obtained by switching the signal of each coordinate. Let  $\{\epsilon_n\}_{n \in \N}$ be a non-increasing sequence of positive numbers such that $\lim_{n \rightarrow \infty} \epsilon_n = \epsilon$ with $0<\epsilon < 1/2$. We define $P$ by
\begin{align*}
&P(1|a^{-1}_{-n}(+\underline{1})) = 1-\epsilon_n = P(-1|-a^{-1}_{-n}(-\underline{1})) \\
&P(1|a^{-1}_{-n}(-\underline{1})) = \epsilon_n = P(-1|-a^{-1}_{-n}(+\underline{1})),
\end{align*}
and put $P(1|\underline{x})=1/2$ for all the remaining pasts $\underline{x}$. 
Clearly $P$ is strongly non-null, attractive, and non-continuous. Also, let
\begin{align*}
&P_+(1|\underline{a}) = \lim_{n \rightarrow \infty} P(1|a^{-1}_{-n}(+\underline{1})) = 1-\epsilon \\
&P_-(-1|\underline{a}) = \lim_{n \rightarrow \infty} P(-1|a^{-1}_{-n}(-\underline{1})) = 1-\epsilon.
\end{align*}
By Lemma 2.3 in \cite{hulse/1991}, the maximum phase ${\bf X}^{(+\underline{1})}$ is consistent with $P_+$ and therefore it is the i.i.d. process with probability $1-\epsilon$ for $1$. Analogously, the minus phase ${\bf X}^{(-\underline{1})}$ is the i.i.d. process with probability $1-\epsilon$ for $-1$. \\

Theorem \ref{theo:attractive} is a direct consequence of Theorems \ref{theo:perfect}, \ref{theo:exponential}, and \ref{theo:unique} below.

\begin{theo} \label{theo:perfect}
Let $P$ be an attractive continuous kernel. If there exists a unique stationary chain compatible with $P$ then there exists a feasible CFTP algorithm using a countable alphabet i.i.d. process that simulates this chain.
\end{theo}

Observe that for  this theorem  we do not require strong non-nullness of the kernel. As an immediate application of Theorem \ref{theo:perfect}, our \emph{Attractive Sampler} given in Section \ref{sec:attractivesampler} perfectly simulates binary auto-regressive and BK processes introduced in Section \ref{sec:def} in their uniqueness regime.
Notice that for the binary autoregressive and BK processes, for any $\eta>0$, we can exhibit kernels having continuity rate $\textrm{var}_{k}=O(1/k^\eta)$ for which the unique compatible stationary chain can be perfectly simulated. In particular, in the BK example $\textrm{var}_{k}$ can be taken so that it converges arbitrarily slowly to $0$.
As a comparison, in the work  of \citet{comets/fernandez/ferrari/2002} the condition $\sum_{k\geq0}\prod_{i=0}^{k-1}(1-\textrm{var}_{k})=+\infty$ is assumed to guarantee that their  CFTP algorithm is feasible, \textit{i.e.}, the stopping time is a.s. finite. This condition does not hold if $\textrm{var}_{k}=O(1/k^\eta)$ with sufficiently small $\eta$. In other words, in the class of regular attractive chains, our perfect simulation algorithm (Attractive Sampler) is optimal. This is particularly clear for the binary auto-regressive model with $\phi(r) = (1+r)/2$, which is also considered in \cite{comets/fernandez/ferrari/2002}. In this case, when $\sum_{n\geq1}\xi_{n} + \gamma <1$, the criterion of \citet{fernandez/maillard/2005} implies uniqueness, and therefore, our Attractive Sampler works whereas the algorithm in \citet{comets/fernandez/ferrari/2002} is not guaranteed to work  in general.\\

The following theorem relates the stopping times of a FP to its concentration property.

\begin{theo} \label{theo:exponential}
If ${\bf X}$ is a FP then the concentration of measure holds at exponential rate for ${\bf X}$.  Explicitly,  using the same notation as in (\ref{eq:checkable}) and (\ref{eq:expoconcentration}), if $\theta_1$ and $\theta_2$ are the stopping times of the FP and $r_1$ and $r_2$ are positive numbers such that 
\begin{equation*}
\P(\{\theta_1 >  r_1\} \cup \{\theta_2 >  r_2\}) \leq \epsilon/(6\|\delta f\|_{\ell_1(\N)}),
\end{equation*} 
then we have 
\begin{equation*}
\P\left(\left|f(X^{n}_{1}) - \E[f(X^n_1)]\right| > \epsilon \right) \leq 4\exp\left\{-\frac{ 2\epsilon^2}{9(1+r_2+\sum_{j=1}^{r_1}\P(\theta_1 \geq j))^{2}\|\delta f \|^2_{\ell_2(\N)} } \right\}.
\end{equation*}
\end{theo}

The above theorem holds for any FP process (we assume neither regularity nor attractiveness) and it is of independent interest. We note that, instead of Theorem \ref{theo:exponential}, the blowing up property of FP proved in \cite{marton/shields/1994} could be used together with our Theorems \ref{theo:perfect} and \ref{theo:unique} to prove the equivalence between uniqueness of compatible chain and the existence of a finitary coding from a countable alphabet i.i.d. process to the compatible chain. This is because in \cite{marton/shields/1994} it is proved that if a process is FP then it satisfies the blowing up property, which implies that the ergodic theorem holds at exponential rate. Nevertheless, we think that Theorem \ref{theo:exponential} gives us more explicit information about the process relating the stopping time to the concentration of measure. For instance, if a process can be sampled using a CFTP algorithm and the stopping time has finite expectation, we have the following useful corollary, which is obtained simply taking $r_1 = \infty$ and $r_2 = 0$.

\begin{coro} \label{prop:cftpconcentration}
Let ${\bf X}$ be a process that can be simulated by a CFTP algorithm with a stopping time $\theta$. If $\E[\theta] < \infty$, then for all $\epsilon > 0$ and all functions $f:A^n \to \R$ we have
\begin{equation} \label{eq:gaussianconcentration}
\P\left(\left|f(X^{n}_{1}) - \E[f(X^n_1)]\right| > \epsilon \right) \leq 4\exp\left\{{-\frac{2\epsilon^2}{9(1+\E[\theta])^2
\|\delta f \|^2_{\ell_2(\N)}}} \right\}.
\end{equation}
\end{coro} 

As an example of application of the above result, if a probability kernel has summable continuity rate $var_k$, we can construct a CFTP algorithm with stopping time $\theta$ such that $\E[\theta] \leq C\sum_{k=1}^\infty var_k$, where $C$ is a numerical constant \cite{comets/fernandez/ferrari/2002}. We refer the reader to \cite{gallo/2011}, \cite{gallo/garcia/2011}, and \cite{desantis/piccioni/2010} to obtain bounds on probability of $\theta$ for different conditions and not necessary regular kernels.

\vspace{0.1cm}

Now we have the last ingredient for the proof of Theorem \ref{theo:attractive}.

\begin{theo} \label{theo:unique}
Let $P$ be a regular kernel and ${\bf X}$ a process compatible with $P$ that satisfies the concentration of measure at exponential rate. Then ${\bf X}$ is the unique stationary process compatible with $P$. 
\end{theo}
Note that, for this result, we do not assume that the kernel is attractive and therefore, Theorem \ref{theo:unique} constitutes an interesting characterization of uniqueness for chains of infinite order. We cannot, in general, relax the strong non-nullness condition in this theorem, because by the example given just after Theorem  \ref{theo:attractive} where the kernel has only one null transition probability, there exists a chain that satisfies the ergodic theorem at exponential rate but it is not the unique chain compatible with the kernel.
  
\begin{figure}[!ht]
\begin{center} 
 \includegraphics[width=12cm]{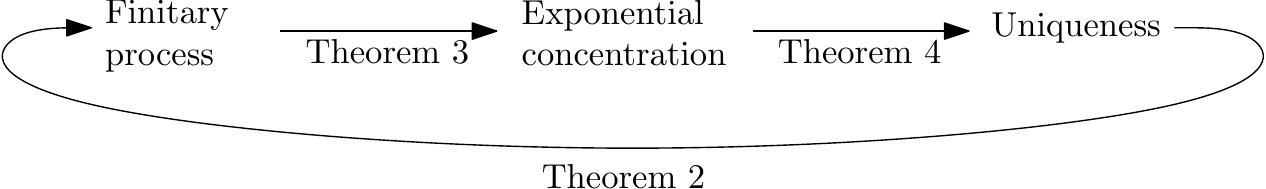}
 \caption{\label{fig:diagram} \small Diagram showing the chain of implications stated by Theorem \ref{theo:attractive}} \label{fig:intmatrix2weak}
 \end{center}
\end{figure}

Now, Theorem  \ref{theo:attractive} follows from the sequence of implications shown in Figure \ref{fig:diagram} that holds when the kernel is attractive and regular.

From the coding point of view, it is natural to ask if Theorem \ref{theo:attractive} can be strengthen to a finitary coding from a finite alphabet i.i.d. process.  This is indeed the case for the original BK example.

\begin{theo} \label{theo:finiteBK}
Let $P$ be the original BK example. Then there exists a unique chain compatible with $P$ if and only if the compatible chain is a finitary coding of a finite alphabet i.i.d. process. 
\end{theo}

\section{The Attractive Sampler}\label{sec:attractivesampler}

Assume that we are given an attractive continuous kernel $P$ for which there exists a unique compatible stationary chain. The alphabet is $A=\{1,\ldots,s\}$. As stated in Theorem \ref{theo:perfect}, there exists a CFTP algorithm that samples from its stationary law. Here we construct one such CFTP algorithm and call it the \emph{Attractive Sampler}.
First, let us consider the kernel $\tilde{P}$ on $A\times A$ introduced in \citep{hulse/1991} and defined by 
\begin{equation}\label{eq:couplingdobom}
\tilde{P}((x_0\ge a,y_0\ge b)|(\underline{x},\underline{y})):=\sum_{i=a}^{s}P(i|\underline{x})\wedge\sum_{i=b}^{s}P(i|\underline{y}),
\end{equation}
for any pair of pasts $\underline{x}$ and $\underline{y}$ in $A^{-\mathbb{N}}$ and any pair of symbols $a$ and $b$ in $A$. Note that, for any $a,b$ in $A$, we have
\begin{align*}
 \tilde{P}((a,b)|(\underline{x},\underline{y})) &= \tilde{P}((x_0=a,y_0=b)|(\underline{x},\underline{y}))\\
& =\tilde{P}((x_0\ge a,y_0\ge b)|(\underline{x},\underline{y})) - \tilde{P}((x_0 \geq a+1,y_0 \geq b)|(\underline{x},\underline{y})) \\
&- \tilde{P}((x_0 \geq a,y_0 \geq b+1)|(\underline{x},\underline{y})) + \tilde{P}((x_0 \geq a+1,y_0 \geq b+1)|(\underline{x},\underline{y})).
\end{align*}
This kernel defines a coupling between the kernel $\{P(a|\underline{x})\}_{a\in A}$ and $\{P(a|\underline{y})\}_{a\in A}$.
To see this, first observe that
\[
\sum_{a=1}^{s}\tilde{P}((a,s)|(\underline{x},\underline{y}))=\tilde{P}((a'\ge 1,b'\ge s)|(\underline{x},\underline{y}))=P(s|\underline{y}),
\]
secondly, observe that 
\begin{align*}
\sum_{a=1}^{s}\tilde{P}((a,s-1)|(\underline{x},\underline{y}))&=\tilde{P}((a'\ge 1,b'\ge s-1)|(\underline{x},\underline{y}))-\tilde{P}((a'\ge 1,b'\ge s)|(\underline{x},\underline{y}))\\
&=\left[P(s-1|\underline{y})+P(s|\underline{y})\right]-P(s|\underline{y})\\
&=P(s-1|\underline{y}).
\end{align*}
Now, continuing the recursion, we show that $\sum_{a=1}^{s}\tilde{P}((a,b)|(\underline{x},\underline{y}))=P(b|\underline{y})$ for any $b\in A$. The same holds for the sum over $b$, that is  $\sum_{b=1}^{s}\tilde{P}((a,b)|(\underline{x},\underline{y}))=P(a|\underline{y})$ for any $a\in A$. This means that $\tilde{P}$ defines a coupling between the chains with respective fixed pasts. 

A straightforward but tedious computation shows that $\tilde{P}$ is indeed continuous and we can use  the result of \citet{kalikow/1990} stating that continuous kernels can be written as a countable mixture of Markov kernels of increasing order. Formally, for $\tilde{P}$, there exists a sequence of non-negative numbers $\{\lambda_{k}\}_{k\ge 0}$ with $\sum_{k\geq 0} \lambda_k = 1$, a sequence of Markov kernels $\{P^{[k]}\}_{k\ge1}$, where $P^{[k]}$ is a $k$-step Markov kernel, and a probability distribution $P^{[0]}$ on $A\times A$ such that
\[
\tilde{P}((a,b)|(\underline{x},\underline{y}))= \lambda_{0}P^{[0]}((a,b))+\sum_{k\ge1}\lambda_{k}P^{[k]}((a,b)|(x_{-k}^{-1},y_{-k}^{-1})).
\]
We will now use this representation of $\tilde{P}$ to define our algorithm. First, we introduce a series of partition of $[0,1]$ that will be used to define the finitary coding.

For any $k \in \N$, $(x^{-1}_{-k},y^{-1}_{-k})\in A^{-k}\times A^{-k}$ and $a$, $b \in \{1, \ldots, s\}$, define
\begin{align*}
r^{[k]}_{a,b}(x^{-1}_{-k}, y^{-1}_{-k}) = &\sum_{c \leq a-1}\sum_{d \leq s}P^{[k]}((x_0 = c,y_0 = d)|(x_{-k}^{-1},y_{-k}^{-1}))\\
&+\sum_{d \leq b}P^{[k]}((x_0 = a,y_0 = d)|(x_{-k}^{-1},y_{-k}^{-1})).
\end{align*}
and
\begin{equation*}
r^{[0]}_{a,b} = \sum_{c \leq a-1}\sum_{d \leq s}P^{[0]}((x_0 = c,y_0 = d))+\sum_{d \leq b}P^{[0]}((x_0 = a,y_0 = d)).
\end{equation*}
 Let $r^{[k]}_{1,0} := 0$ and, for $a \in\{2,\ldots,s\}$, $r^{[k]}_{a,0} := r^{[k]}_{a-1,s}$. Notice that 
 $$r^{[k]}_{a,b}(x^{-1}_{-k}, y^{-1}_{-k}) - r^{[k]}_{a,b-1}(x^{-1}_{-k}, y^{-1}_{-k}) = P^{[k]}((x_0 = a,y_0 = b)|(x_{-k}^{-1},y_{-k}^{-1})).$$
 For any integer $k \geq 0$, we set
\begin{equation}\label{eq:barrinhas}
\mathcal{R}^{[k]} = \left\{ r^{[k]}_{a,b}(x^{-1}_{-k}, y^{-1}_{-k}) \in [0,1]:  (a,b) \in A^2, (x_{-k}^{-1},y_{-k}^{-1}) \in A^{-k}\times A^{-k} \right\}.
\end{equation}
Denote the elements of $\mathcal{R}^{[k]}$ by $ t^{[k]}_{1} < t^{[k]}_{2} < \ldots < t^{[k]}_{|\mathcal{R}^{[k]}|} = 1$ and define $t^{[k]}_{0} := 0$. Now, for all $k \in \N$, define the following set of intervals
\begin{equation*}
\mathcal{I}^{[k]} = \left\{\left[\sum_{i=0}^{k-1}\lambda_i+\lambda_kt^{[k]}_{j-1},\;\; \sum_{i=0}^{k-1}\lambda_i+ \lambda_kt^{[k]}_{j}\right[:  j \in \left\{1, \ldots, |\mathcal{R}^{[k]}|\right\}  \right\}
\end{equation*}
and
\begin{equation*}
\mathcal{I}^{[0]} = \left\{\left[\lambda_0t^{[0]}_{j-1},\;\;  \lambda_0t^{[0]}_{j}\right[:  j \in \left\{1, \ldots, |\mathcal{R}^{[0]}|\right\}  \right\}
\end{equation*}
Finally, let $\mathcal{I} = \bigcup_{k\geq 0} \mathcal{I}^{[k]}$. Observe that  the set $\mathcal{I}$ has a countable number of elements, therefore, we can identify its elements by $I_j, j \in \N$. We define the \emph{update} function $F:A^{-\mathbb{N}}\times A^{-\N}\times \N\rightarrow A\times A$ by 
\begin{align*}\label{eq:updatesimples}
&F((\underline{x}, \underline{y}),j) \\
&= \sum_{(a,b) \in A^2}(a,b){\bf 1}\left\{ I_j \subset \bigcup_{k \in \N}\left[r^{[k]}_{a,b-1}(x^{-1}_{-k}, y^{-1}_{-k}),r^{[k]}_{a,b}(x^{-1}_{-k},y^{-1}_{-k}) \right[   \bigcup \left[r^{[0]}_{a,b-1},r^{[0]}_{a,b}\right[\right\}.
\end{align*}

Now, let $\mathbf{L} = \{L_j\}_{j \in \Z}$ be an i.i.d. process with values on $\N$ such that 
\begin{equation}\label{eq:Lacopla}
\P\left(L_0 = j \right)= \sup I_j - \inf I_j.
\end{equation}
We observe by construction that 
\begin{equation*}
\P(F((\underline{x}, \underline{y}),L_0) = (a,b)) = \tilde{P}((a,b)|(x,y)),
\end{equation*}
which justifies the name update function for $F$.

Let the concatenation of  pairs of symbols  be understood coordinatewise, \textit{i.e.}, $(a,b)(c,d) = (ac, bd)$ whenever $(ac, bd)$ is well defined. Using this notation, we define, for any $-\infty<k\leq l <+\infty$, the successive iterations of $F$ as
\begin{equation*}
F_{[k,l]}((\underline{x}, \underline{y}),L_{k}^{l})=F\left(F_{[k,l-1]}((\underline{x}, \underline{y}),L_{k}^{l-1}) \ldots F_{[k,k]}((\underline{x}, \underline{y}),L_{k})(\underline{x}, \underline{y}),L_{l}\right),
\end{equation*}
where  $F_{[k,k]}((\underline{x},\underline{y}),L_{k}):=F((\underline{x},\underline{y}),L_{k})$. Notice that for any $(\underline{x}, \underline{y})$,
we can obtain a coupling $\{(X^{\underline{x}}_{i}, X^{\underline{y}}_{i})\}_{i\geq0}$ by observing that
$$F_{[0,i]}((\underline{x}, \underline{y}),L_{0}^{i}) \eqd (X^{\underline{x}}_{i}, X^{\underline{y}}_{i}).$$ 
Furthermore, we define for any $i\in\mathbb{Z}$, the random variable
\begin{equation*}
\theta[i]:=\min\{j\geq 0:F_{[i-j,i]}((\underline{a}, \underline{b}),L_{i-j}^{i}) \in \{(1,1), \ldots, (s,s)\}\,\,\,\textrm{for all}\,\,(\underline{a},\underline{b}) \in A^{-2\N}\}.
\end{equation*}
An important observation is that in the particular case of attractive chains, 
\begin{equation*}
\theta:=\theta[0]=\min\{j\geq 0:F_{[-j,0]}((\underline{1}, \underline{s}),L_{-j}^{0}) \in \{(1,1), \ldots, (s,s)\} \}.
\end{equation*}
Finally, define the value of the  \emph{coding function} $\Phi$ at time $i\in\mathbb{Z}$ by
\begin{equation*}\label{eq:Phi3}
[\Phi({\bf L})]_{i}=F_{[-\theta[i],i]}((\underline{1}, \underline{s}),L_{-\theta[i]}^{i}).
\end{equation*}
The algorithm Attractive Sampler is defined by the pseudocode below. 
\begin{algorithm}[h] 
\caption*{Algorithm Attractive Sampler} 
\begin{algorithmic}[1]
\STATE {\it Input:} $F$;\,\,{\it Output:} $\theta[0]$, $[\Phi({\bf L})]_{0}$
\vspace{0.1cm}
\STATE Sample $L_{0}$ with distribution $\mathbb{P}$
\STATE $i \leftarrow 0$, $\theta[0] \leftarrow 0$, $[\Phi({\bf L})]_{0}\leftarrow 0$
\WHILE{$F_{[-i,0]}((\underline{1}, \underline{s}),L_{-i}^{0})\notin \{(1,1), \ldots, (s,s)\}$}
\STATE $i \leftarrow i+1$\\
\STATE Sample $L_{-i}$ with distribution $\mathbb{P}$\\
\ENDWHILE
\STATE $\theta[0]\leftarrow i$
\STATE $[\Phi({\bf L})]_{0}\leftarrow  F_{[-i,0]}((\underline{1}, \underline{s}),L_{-i}^{0})$
\RETURN $\theta[0]$, $[\Phi({\bf L})]_{0}$.

\end{algorithmic}
\end{algorithm}

Observe that definitions above of $\theta$ and $\Phi$ satisfy the requirements of a stationary coding from ${\bf L}$. In particular, it can be shown in a standard way (see for example \cite{propp/wilson/1996} for the Markovian case, or \cite{comets/fernandez/ferrari/2002} for chains of infinite order) that if $\theta$ is $\mathbb{P}$-a.s. finite (that is, the CFTP algorithm is feasible), then the coding is finitary and the sample $[\Phi({\bf L})]_0$ is constructed according to the unique stationary measure compatible with $P$.  The compatibility and the stationarity follow from the construction and the property of the update function. The uniqueness follows from the loss of memory the chain inherits because of the existence of almost 
surely finite stopping time $\theta[i]$ for any $i\in\mathbb{Z}$.

 In Section \ref{sec:prooftheoperfect} (proof of Theorem \ref{theo:perfect}), we will prove that for attractive continuous $P$, the Attractive Sampler is feasible if $P$ has a unique compatible chain.

\section{Proof of the results} \label{sec:proof}

From Section  \ref{sec:results} it is clear that Theorem \ref{theo:attractive} follows directly from Theorems \ref{theo:perfect},  \ref{theo:exponential}, and \ref{theo:unique}.

\subsection{Proof of Theorem \ref{theo:perfect}}\label{sec:prooftheoperfect}
First we need the following lemma proved for the regular and attractive kernels by \citet{hulse/1991}. Here, we drop the unnecessary non-nullness condition.
\begin{lemma} \label{lemma:uniqueness}
Let $P$ be attractive and continuous and consider the update function $F$ defined in Section \ref{sec:attractivesampler}. If there exists a unique chain compatible with $P$, then for all $i \in A = \{2, \ldots, s\}$,  
\begin{equation} \label{conv}
\lim_{n \rightarrow\infty}\mathbb{P}\left(F_{[0,n]}((\underline{1},\underline{s}),L_{0}^{n})\geq (i, 1)\right)-\mathbb{P}\left(F_{[0,n]}((\underline{1},\underline{s}),L_{0}^{n})\geq (1, i)\right) = 0.
\end{equation}
\end{lemma}

\begin{obs} 
Remember that $F_{[0,n]}((\underline{1}, \underline{s}),L_{0}^{n}) \eqd (X^{\underline{1}}_{n}, X^{\underline{s}}_{n})$ and, therefore, we have that  $\mathbb{P}\left(F_{[0,n]}((\underline{1},\underline{s}),L_{0}^{n})\geq (i, 1)\right)$, $i\in A$,
gives the law of $X^{\underline{1}}_{n}$.
\end{obs}

\proc{Proof.}
 For $a,b,c,d \in A$, we write $(a,b) \geq (c,d)$ if $a \geq c$ and $b \geq d$. 
 
 Because $P$ is attractive, for all $\l_1, \ldots, l_k, k \in \N$ and $a_1, \ldots, a_k \in A$, we have that 
 $$\P(F_{[0,n+l_1]}((\underline{1}, \underline{s}),L_{0}^{n+l_1}) \geq (1, a_1),\ldots, F_{[0,n+l_k]}((\underline{1}, \underline{s}),L_{0}^{n+l_k}) \geq (1,a_k) )$$ and 
  \begin{equation*} \label{eq:minimumphase}
 \P(F_{[0,n+l_1]}((\underline{1}, \underline{s}),L_{0}^{n+l_1}) \geq (a_1,1),\ldots, F_{[0,n+l_k]}((\underline{1}, \underline{s}),L_{0}^{n+l_k}) \geq (a_k,1) )
 \end{equation*}
 are respectively non-increasing and non-decreasing in $n \in \N$. Therefore, both sequences are convergent in $n \in \N$ and their limits when $n$ diverges
 define  stationary chains. By construction, for all $b\in A$ and $\underline{a} \in A^{-\N}$ the respective chains are compatible with the kernels $P^{\max}$  given by $P^{\max}(b| \underline{a}) = \lim_{n \rightarrow \infty} P(b| a^{-1}_{-n}\underline{s})$ and $P^{\min} (b| \underline{a}) = \lim_{n \rightarrow \infty} P(b| a^{-1}_{-n}\underline{1})$.

If the kernel $P$ is continuous, $P^{\max} = P^{\min} = P$, and therefore both chains are compatible with $P$. This implies that if there exists only one chain compatible with $P$, then we have, for $i=2,\ldots,s$, 
\begin{equation*}
\lim_{n \rightarrow\infty}\mathbb{P}\left(F_{[0,n]}((\underline{1},\underline{s}),L_{0}^{n})\geq (i, 1)\right)-\mathbb{P}\left(F_{[0,n]}((\underline{1},\underline{s}),L_{0}^{n})\geq (1, i)\right) = 0,
\end{equation*} 
as we wanted.

\ep
\medbreak

We now prove that, for attractive regular chains, convergence (\ref{conv})  implies that the unique stationary chain compatible with $P$ can be sampled by
the algorithm Attractive Sampler. In other words, we need to prove that uniqueness of compatible chain implies that $\theta[0]$ is $\mathbb{P}$-a.s. finite. We have
\[
\mathbb{P}(\theta[0]>n)=\mathbb{P}\left(\cap_{j=-n}^{0}\left\{F_{[j,0]}((\underline{1}, \underline{s}),L_{j}^{0}) \notin \{(1,1), \ldots, (s,s)\}\right \}\right),
\]
which yields, using first the attractiveness and then the translation invariance of ${\bf L}$
\begin{align*}
\mathbb{P}(\theta[0]>n)&= \mathbb{P}\left(F_{[-n,0]}((\underline{1}, \underline{s}),L_{-n}^{0}) \notin \{(1,1), \ldots, (s,s)\}\right)\\
&=\mathbb{P}\left(F_{[0,n]}((\underline{1}, \underline{s}),L_{0}^{n}) \notin \{(1,1), \ldots, (s,s)\}\right).
\end{align*}
Thus, 
we want to prove  that
\[
\lim_{n\rightarrow+\infty}\mathbb{P}\left(F_{[0,n]}((\underline{1}, \underline{s}),L_{0}^{n}) \notin \{(1,1), \ldots, (s,s)\}\right)=0.
\]
Due to the attractiveness, our coupling guarantees that for $a\in A$
\[
\mathbb{P}\left(F_{[0,n]}((\underline{1}, \underline{s}),L_{0}^{n})\geq (a,a) \right)=\mathbb{P}\left(F_{[0,n]}((\underline{1}, \underline{s}),L_{0}^{n})\geq (a,1) \right).\]
Taking the limit and using (\ref{conv}), we have for any $a\in A$ 
\begin{align*}
\lim_{n\rightarrow+\infty}\mathbb{P}\left(F_{[0,n]}((\underline{1}, \underline{s}),L_{0}^{n})\geq (a,a)\right)&=\lim_{n\rightarrow+\infty} \mathbb{P}\left(F_{[0,n]}((\underline{1}, \underline{s}),L_{0}^{n}) \geq (a, 1) \right)\\
&=\lim_{n\rightarrow+\infty} \mathbb{P}\left(F_{[0,n]}((\underline{1}, \underline{s}),L_{0}^{n}) \geq (1, a) \right).
\end{align*}
Now, for any $a \in A$, let $\Gamma(a) = \{(i,j) \in A^2: i \geq a \;\;\text{and}\;\; j < a \}$. From the last equation, we have that
\begin{equation*}
\lim_{n\rightarrow+\infty}\mathbb{P}\left(F_{[0,n]}((\underline{1}, \underline{s}),L_{0}^{n}) \in \Gamma(a) \right) = 0,
\end{equation*}
and this implies
\begin{equation*}
\lim_{n\rightarrow+\infty}\mathbb{P}\left(F_{[0,n]}((\underline{1}, \underline{s}),L_{0}^{n}) \notin \{(1,1), \ldots, (s,s)\}\right)=0,
\end{equation*}
which concludes the proof.

\subsection{Proof of Theorem \ref{theo:exponential}}

Let ${\bf X}$ be any process with alphabet $A$. For all $\underline{x} \in A^{-\N}$, let ${\bf X}^{\underline{x}}$ be the process with fixed past $\underline{x}$. Let $\underline{y}, \underline{z} \in A^{-\N}$ and $({\bf X}^{\underline{y}}, {\bf X}^{\underline{z}})$ be a coupling between the process ${\bf X}^{\underline{y}}$ and ${\bf X}^{\underline{z}}$. The following lemma, which we state without proof, is a direct consequence of Theorem 1 of \cite{chazottes/collet/kulske/regig/2000}.
\begin{lemma}[{\citet{chazottes/collet/kulske/regig/2000}}] \label{theo:chazottes}
Let $({\bf X}^{\underline{y}}, {\bf X}^{\underline{z}})$ be couplings for each  pair  $\underline{y}, \underline{z} \in A^{-\N}$. If $\sup_{\underline{y}, \underline{z}}\sum_{j =1}^\infty \P(X^{\underline{y}}_j \neq X^{\underline{z}}_j)\leq \Delta < \infty$, then for all integer $n \geq 1$, functions $f:A^n \to \R$, and $ \epsilon > 0$ we have
\begin{equation*} 
\P(|f(X^n_1)- \E[f(X^n_1)]| > \epsilon) \leq 2\exp\left\{-\frac{2\epsilon^2}{(1+\Delta)^2\|\delta f\|_{\ell_2(\N)}^2}\right\}.
\end{equation*}
\end{lemma}

Assume that the process ${\bf X}$ is a finitary coding of a sequence  ${\bf U}$, where $U_{i}\in\mathcal{U}$ for any $i\in\mathbb{Z}$. Let $\theta_1, \theta_2$, and $\Phi$ be the quantities involved in the finitary coding as defined generically in Section \ref{sec:def} (and not necessarily as in Section \ref{sec:attractivesampler}, which is specific for algorithm Attractive Sampler). 
Let $\epsilon > 0$.  Take two positive numbers $r_1$ and $r_2$ such that 
\begin{equation*}
\P(\{\theta_1 >  r_1\} \cup \{\theta_2 >  r_2\}) \leq \epsilon/(6\|\delta f\|_{\ell_1(\N)}).
\end{equation*}

Now, we want to show that we can approximate the FP by two $(r_1+r_2)$-dependent processes, for which we have good control of the concentration of measure property.

For $j \in \Z$, let ${\bf V}^{[j]}$ be a family of i.i.d. processes with values in $\mathcal{U}$ independent of each other and of ${\bf U}$. We define a process ${\bf Y}$ as
\begin{equation*}
Y_j = [\Phi(V^{[j],\infty}_{j+r_2}U^{j+r_2}_{j-r_1}V^{[j], j-r_1}_{-\infty})]_j\,\,\,\,,\,\,\,\,\,\textrm{for any}\,\, j\in\Z
\end{equation*}
where, for any $i$ and $j$ in $\Z$, we use the notation $V^{[j],\infty}_{i}$ for the sequence $\ldots V^{[j]}_{i+1}\,V^{[j]}_{i}$ and the notation $V^{[j],i}_{-\infty}$ for the sequence $V^{[j]}_{i-1}V^{[j]}_{i-2} \ldots$. 
Clearly, ${\bf Y}$ is stationary and if $\theta_1\leq r_1$ and $\theta_2\leq r_2$, then $Y_0 = X_0$. Moreover, ${\bf Y}$ is a $(r_1+r_2)$-dependent process  \textit{i.e.}, for  all $l,m > 1$ and ${\bf y} \in A^{\Z}$
\begin{equation*}
\P(Y^{l}_{1} = y^{l}_{1}, Y^{l+r_1+r_2+m}_{l+r_1+r_2+1}=y^{l+r_1+r_2+m}_{l+r_1+r_2+1}) = \P(Y^{l}_{1} = y^{l}_{1}) \P(Y^{l+r_1+r_2+m}_{l+r_1+r_2+1}=y^{l+r_1+r_2+m}_{l+r_1+r_2+1}).
\end{equation*}
We will now use the following equality
\begin{equation} \label{eq:dec}
f(X^{n}_{1} )- \E[f(X^{n}_{1} )] = f(X^{n}_{1} )  - f(Y^{n}_{1} ) - \E[f(X^{n}_{1} )-f(Y^{n}_{1} )] +  f(Y^{n}_{1} ) - \E[f(Y^{n}_{1} )].
\end{equation}
Consider the events $O_1^j = \{\theta_1(T_{\mathcal{U}}^j{\bf U}) > r_1\}$ and $O_2^j = \{\theta_2(T_{\mathcal{U}}^j{\bf U}) > r_2\}$. From the definition of $\delta f$, we have that 
\begin{equation} \label{eq:decineq1}
\left |f(X^{n}_{1} )  - f(Y^{n}_{1} )\right| \leq \sum_{j=1}^n{\bf 1}\{X_{j} \neq Y_{j} \}\delta_jf \leq \sum_{j=1}^n{\bf 1}\{O_1^j \cup O_2^j \}\delta_jf
\end{equation}
 and 
\begin{equation}\label{eq:decineq2}
\left| \E[f(X^{k}_{1} )-f(Y^{k}_{1} )]  \right | \leq  \E\left[\left|f(X^{k}_{1} )-f(Y^{k}_{1} ) \right | \right]\leq \P(O_1^0 \cup O_2^0 )\|\delta f\|_{\ell_1(\N)} \leq \epsilon/6.
\end{equation}

Observe that we have
\begin{align*}
&\P\left(|f(X^{n}_{1} )- \E[f(X^{n}_{1} )] | > \epsilon\right) \\ 
&\leq \P\left(|f(X^{n}_{1} )  - f(Y^{n}_{1} ) - \E[f(X^{n}_{1} )-f(Y^{n}_{1} )]|  > 2\epsilon/3\right) \notag \\
& + \P\left(| f(Y^{n}_{1} ) - \E[f(Y^{n}_{1} )]| > \epsilon/3\right).
\end{align*}

Therefore, collecting (\ref{eq:dec}), (\ref{eq:decineq1}), (\ref{eq:decineq2}), we have 
\begin{align}
&\P\left(|f(X^{n}_{1} )- \E[f(X^{n}_{1} )] | > \epsilon\right) \notag \\ 
&\leq \P\left(\left|\sum_{j=1}^n{\bf 1}\{O_1^j \cup O_2^j \}\delta_jf -  \P(O_1^0 \cup O_2^0)\|\delta f\|_{\ell_1(\N)}  \right| > \epsilon/3\right) \notag \\
& + \P\left(| f(Y^{n}_{1} ) - \E[f(Y^{n}_{1} )]| > \epsilon/3\right). \label{eq:conc}
\end{align}
We will use Lemma \ref{theo:chazottes} to obtain upper bounds for the two terms  of the  right hand side of (\ref{eq:conc}). 

Let us begin with the second term, and we will use the fact that ${\bf Y}$ is a $r_1+r_2$-dependent process. 
First we define a coupling $({\bf \tilde{Y}}, {\bf \hat{Y}})$, where ${\bf \tilde{Y}}$ and ${\bf \hat{Y}}$ are copies of ${\bf Y}$, by
$$(\tilde{Y}_j, \hat{Y}_j) := \left([\Phi(V^{[j],\infty}_{j+r_2}U^{j+r_2}_{j-r_1}V^{[j], j-r_1}_{-\infty})]_j, [\Phi(V'^{[j],\infty}_{j+r_2}U'^{j+r_2}_{j-r_1}V'^{[j], j-r_1}_{-\infty})]_j\right), \;\;\text{for all}\;\;\;  j \in \Z$$
 where ${\bf V}^{[j]}, {\bf V'}^{[j]}$ and ${\bf U}, {\bf U'}$ are i.i.d. processes satisfying the following properties. For $j > 0$, ${\bf V}^{[j]} = {\bf V'}^{[j]}$. For $j \leq 0$, ${\bf V}^{[j]}$ is independent of  ${\bf V'}^{[j]}$, and both are independent of the rest. The processes ${\bf U}, {\bf U'}$ are independent of ${\bf V}^{[j]}, {\bf V'}^{[j]}$ for all $j \in \Z$. Also $\{U_j\}_{j\leq r_2}$ and $\{U'_j\}_{j\leq r_2}$ are independent and for $j > r_2$, $U_j = U'_j$.
From the construction, for all $\underline{y},\underline{z} \in A^{-\N}$ and $j > r_1+r_2$,
 \begin{equation*}
 \P(\tilde{Y}_j =  \hat{Y}_j | \tilde{Y}^0_{-\infty} = \underline{y},  \hat{Y}^0_{-\infty} = \underline{z}) = \P(\tilde{Y}_j = \hat{Y}_j)=1.
 \end{equation*}
 
Therefore, using the above coupling and the fact that ${\bf Y}$ is a $(r_1+r_2)$-dependent process, we have $\Delta = r_2+\sum_{j=1}^{r_1}\P(\theta_1\geq j)$ in Lemma \ref{theo:chazottes}. Hence, we have
\begin{equation} \label{eq:firstconc1}
\P\left(| f(Y^{n}_{1} ) - \E[f(Y^{n}_{1} )]| > \epsilon/3\right) \leq 2\exp\left\{-\frac{g(\epsilon, \|\delta f\|_{\ell_1(\N)})}{\|\delta f\|_{\ell_2(\N)}^2} \right\}
\end{equation}
where $g(\epsilon, \|\delta f\|_{\ell_1(\N)}) = \frac{2}{9}\epsilon^2(1+r_2+\sum_{j=1}^{r_1}\P(\theta_1 \geq j))^{-2}$. Note that $r_1$ and $r_2$ depend on $\epsilon$ and $\|\delta f\|_{\ell_1(\N)} $.

For the first term of the  right hand side of (\ref{eq:conc}),  let $Z_j = {\bf 1}\{O_1^j \cup O_2^j \}$. Observe that the process ${\bf Z} = \{Z_j\}_{j \in \Z}$ is stationary. Also because the event $O_1^j \cup O_2^j $ is $\mathcal{F}(U^{r_2+j}_{-r_1+j})$ measurable, the process ${\bf Z}$ is a $(r_1+r_2)$-dependent process. We define a coupling $(\tilde{\bf Z},\hat{\bf Z})$, where ${\bf \tilde{Z}}$ and ${\bf \hat{Z}}$ are copies of ${\bf Z}$, by 
\begin{align*}
&( \tilde{Z}_j, \hat{Z}_j) \\
&= \left({\bf 1}\left\{\{\theta_1(T_{\mathcal{U}}^j{\bf U}) \geq r_1\} \cup \{\theta_2(T_{\mathcal{U}}^j{\bf U}) \geq r_2\} \right\}, {\bf 1}\left\{\{\theta_1(T_{\mathcal{U}}^j{\bf U'}) \geq r_1\} \cup \{\theta_2(T_{\mathcal{U}}^j{\bf U'}) \geq r_2\} \right\}\right)
\end{align*}
where ${\bf U}$ and ${\bf U}'$ are i.i.d. processes with $U_j=U'_j$ for $j > r_2$ and $\{U_j\}_{j\leq r_2}$ and $\{U'_j\}_{j\leq r_2}$ are independent.  Therefore, in the same way as (\ref{eq:firstconc1}), by Lemma \ref{theo:chazottes} we have that
\begin{equation*} 
\P\left(\left|\sum_{j=1}^n{\bf 1}\{O_1^j \cup O_2^j \}\delta_jf -  \P(O_1^0 \cup O_2^0)\|\delta f\|_{\ell_1(\N)}  \right| > \frac{\epsilon}{3}\right)\leq 2\exp \left\{-\frac{g(\epsilon, \|\delta f\|_{\ell_1(\N)})}{\|\delta f\|_{\ell_2(\N)}^2}\right\} ,
\end{equation*}
where $g$ was defined in (\ref{eq:firstconc1}). Finally, we have
\begin{equation*}
\P\left(|f(X^{n}_{1} )- \E[f(X^{n}_{1} )] | > \epsilon\right) \leq 4\exp\left\{-\frac{g(\epsilon, \|\delta f\|_{\ell_1(\N)})}{\|\delta f\|^2_{\ell_2(\N)}}\right\},
\end{equation*}
which proves the theorem.

\subsection{Proof of Theorem \ref{theo:unique}}

We say that a stationary process ${\bf X}$ has the \emph{positive divergence property} if
\begin{equation*}
\liminf_{n\rightarrow+\infty}\frac{1}{n+1}\sum_{a \in A^{n+1}}\P(Y^{0}_{-n}=a^{0}_{-n})\log{\frac{\P(Y^{0}_{-n}=a^{0}_{-n})}{\P(X^{0}_{-n}=a^{0}_{-n})}} > 0
\end{equation*}
for any ergodic process ${\bf Y}$ different of ${\bf X}$. 

The proof of Theorem \ref{theo:unique} is based on the following lemmas. The first result is from Theorem 1 in \citet{marton/shields/1994}.
\begin{lemma}[\citet{marton/shields/1994}] \label{lemma:marton}
Let ${\bf X}$ be a stationary process with ergodic theorem holding at exponential rate. Then ${\bf X}$ has the positive divergence property.
\end{lemma}

Now, we will prove that two stationary chains compatible with same regular probability kernel cannot be distinguished by the divergence rate.
\begin{lemma} \label{lemma:relativeentropy}
Let ${\bf X}$ and ${\bf Y}$ be two stationary processes compatible with a continuous kernel $P$. Let $\inf_{\underline{a} \in A^{\Z}}P(a_0|a^{-1}_{-\infty}) = \delta > 0$. Then the relative entropy rate
\begin{equation*}
\lim_{n \rightarrow \infty}\frac{1}{n+1}\sum_{a \in A^{n+1}}\P(X^{0}_{-n}=a^{0}_{-n})\log{\frac{\P(X^{0}_{-n}=a^{0}_{-n})}{\P(Y^{0}_{-n}=a^{0}_{-n})}}
\end{equation*}
exists and is 0.
\end{lemma}

\proc{Proof.}
Let ${\bf Z}$ represent ${\bf X}$ or ${\bf Y}$. Define, for $k \in \N$, 
\begin{equation*}
H_{\bf X}(Z^{0}_{-k}) = -\sum_{a \in  A^{k+1}} \P(X^{0}_{-k}=a^{0}_{-k})\log{\P(Z^{0}_{-k} = a^{0}_{-k})}.
\end{equation*}
Now, we can rewrite the relative entropy rate as
\begin{equation*}
\lim_{n \rightarrow \infty}\frac{1}{n+1}\left\{H_{\bf X}(Y^{0}_{-{n}})-H_{\bf X}(X^{0}_{-{n}})\right\}.
\end{equation*}
Define also, for $k \in  \N$,
\begin{equation*}
H_{\bf X}(Z_0|Z^{-1}_{-k}) = -\sum_{a \in A^{k+1}} \P(X^{0}_{-k}=a^{0}_{-k})\log{\P(Z_0 = a_0|Z^{-1}_{-k} = a^{-1}_{-k})}.
\end{equation*}
By the chain rule and the stationarity of the processes, we have
\begin{equation*}
H_{\bf X}(Z^{0}_{-k}) = H_{\bf X}(Z_0)+\sum_{k=1}^{n} H_{\bf X}(Z_0|Z^{-1}_{-k}).
\end{equation*}
Therefore, we have that the relative entropy rate is a difference between two Ces\`aro sums. To prove that the relative entropy exists, it is enough to show that the limit
\begin{equation*}
\lim_{n \rightarrow \infty}H_{\bf X}(Z_0|Z^{-1}_{-n})
\end{equation*}
exists.
To see that $H_{\bf X}(Z_0|Z^{-1}_{-n})$ converges, let $\mu_{\bf Z}$ be the measure associated with ${\bf Z}$, we have
\begin{equation*}
H_{\bf X}(Z_0|Z^{-1}_{-n}) = -\E_{\bf X}(\log(\mu_{\bf Z}(X_{0}|X^{-1}_{-n}))).
\end{equation*}
By assumption, for all $a \in A^{\Z}$
\begin{equation*}
-\log(\mu_{\bf Z}(a_{0}|a^{-1}_{-n})) \leq -\log \delta,
\end{equation*}
therefore, by the dominated convergence theorem
\begin{equation*}
\lim_{n \rightarrow \infty}-\E_{\bf X}(\log(\mu_{\bf Z}(X_{0}|X^{-1}_{-n}))) = -\E_{\bf X}(\lim_{n \rightarrow \infty}\log(\mu_Z(X_{0}|X^{-1}_{-n}))).
\end{equation*}
By continuity of $P$ we have that, for all $a \in A^{\Z}$
\begin{equation*}
\lim_{n \rightarrow \infty}\log(\mu_{\bf Z}(a_{0}|a^{-1}_{-n})) = \log P(Y_0=a_0|Y^{-1}_{-\infty} = a^{-1}_{-\infty}).
\end{equation*}
Hence,
\begin{equation*}
\lim_{n \rightarrow \infty} H_{\bf X}(Z_0|Z^{-1}_{-n}) = -\E_{\bf X}(\log(P(X_{0}|X^{-1}_{-\infty}))),
\end{equation*}
which concludes the proof. 
\ep
\medbreak

\proc{Proof of Theorem \ref{theo:unique}.}
If ${\bf X}$ has the concentration of measure holding at exponential rate, then we have that it satisfies the ergodic theorem at exponential rate (see equation (\ref{eq:exfreq})). By Lemma \ref{lemma:marton}, ${\bf X}$ has the positive divergence property. By Lemma \ref{lemma:relativeentropy}, if ${\bf X}$ has the positive divergence property, there is no other ergodic process compatible with $P$. Therefore, we conclude that ${\bf X}$ is the unique stationary process compatible with $P$.
\ep
\medbreak

\subsection{Proof of Theorem \ref{theo:finiteBK}}
The kernel of the original BK example (see Section \ref{sec:def}) is defined through
\[
P(+1|\underline{x})=\sum_{k\geq1}\lambda_{k}\phi\left(\frac{1}{m_{k}}\sum_{i=1}^{m_{k}}x_{-i}\right).
\]
with $\phi(t)=(1-\epsilon){\bf 1}\{t>0\}+\epsilon{\bf 1}\{t<0\}$ for some $\epsilon\in(0,1/2)$, an increasing sequence of odd positive integers $\{m_j\}_{j \in \N}$, and $\lambda_k = cr^{k}$ for some $r \in (2/3,1)$ and $c=(1-r)/r$.

In other words, the BK example is given under the form of a countable mixture of Markov kernels $P^{[m_k]}$, $k\ge1$, of lacunary ranges. 
When uniqueness holds, we denote by ${\bf X}^{BK}$ the stationary chain compatible with $P$.

The proof of Theorem \ref{theo:finiteBK} is based on Lemma \ref{lemma:1bis} below, which implies that the ${\bf X}^{BK}$ is a finitary coding of a \emph{finite entropy} i.i.d. process with countable alphabet. Then, we can use Theorem 8 in \cite{rudolph/1982} which states that countable state mixing Markov process having an exponentially decaying return times state is finitarily Bernoulli. This implies, in particular, that  i.i.d. chains with finite entropy are finitary codings of i.i.d. chains on finite alphabet.  We  conclude the proof of the theorem observing that  if a process ${\bf X}$ is a finitary coding of a process ${\bf Y}$, which is itself a finitary coding of a process ${\bf Z}$, then ${\bf X}$ is also a finitary coding of ${\bf Z}$.

\begin{lemma}\label{lemma:1bis}
The sequence ${\bf L}$ used in the attractive sampler for ${\bf X}^{BK}$ can be chosen to have finite entropy.
\end{lemma}

\proc{Proof of Lemma \ref{lemma:1bis}.}
The sequence ${\bf L}$ used in the attractive sampler (see Section \ref{sec:attractivesampler}) derives from the choice of a countable decomposition of the (maximal) coupling kernel $\tilde{P}$. In the present case, the special form of the BK example will  naturally yield the decomposition, as it is itself defined through such decomposition. 

Let us first rewrite the kernel $P$ in the following way
\begin{equation}\label{eq:decompoBK2}
P(+1|\underline{x})=\epsilon+\sum_{k\geq1}(1-2\epsilon)\lambda_{k}{\bf 1}\{\textrm{maj}(x_{-m_k}^{-1})=1\},
\end{equation}
where $\textrm{maj}(x_{-j}^{-1}) = {\bf 1}\{\sum_{i=1}^jx_{-i} \geq 0\} - {\bf 1}\{\sum_{i=1}^jx_{-i} < 0\}$.
Representation (\ref{eq:decompoBK2}) motivates the following way to write the coupling kernel:
\begin{equation*}
\tilde{P}((a,b)|(\underline{x},\underline{y}))=\tilde\lambda_{0}\tilde{P}^{[0]}((a,b))+\sum_{k\ge1}\tilde\lambda_{k}\tilde{P}^{[k]}((a,b)|(x_{-k}^{-1},y_{-k}^{-1})),
\end{equation*}
where 
\begin{itemize}
\item $\tilde{\lambda}_{0}=2\epsilon$,  $\tilde{\lambda}_{m_{k}}=(1-2\epsilon)\lambda_{k}$ for $k\ge1$ and $\tilde\lambda_{i}=0$ for any $i\notin\{0, m_{1},m_{2},\ldots\}$,
\item $\tilde P^{[0]}((a,b)):=\frac{1}{2}{\bf 1}\{a=b\}$,
\item and $\tilde{P}^{[k]}((a,b)|(x_{-k}^{-1},y_{-k}^{-1}))={\bf 1}\{a=\textrm{maj}(x_{-k}^{-1})\;\text{and}\;b=\textrm{maj}(y_{-k}^{-1})\}$.
\end{itemize}
We observe that this kernel satisfies \eqref{eq:couplingdobom}, and that the Markovian kernels of the decomposition are deterministic with respect to the pasts (only $0$'s and $1$'s in the transition probabilities). It follows that, for any $k\ge1$, the set $\mathcal{R}^{[m_k]}$ (see \eqref{eq:barrinhas}) consists of only two values, $\mathcal{R}^{[m_k]}=\{0,1\}$ and 
$$\mathcal{I}^{[m_k]} = \left\{ \left[\sum_{j=1}^{k-1}\tilde \lambda_{m_j},\; \sum_{j=1}^{k}\tilde \lambda_{m_j}\right[\right\}.$$
We also observe that 
$$\mathcal{I}^{[0]} = \left\{ \left[0,\; \epsilon\right[; \left[\epsilon,\; 2\epsilon\right[\right\}.$$
Therefore, the entropy of the i.i.d. process $\mathbf{L}$ is 
$\epsilon\log \epsilon + \sum_{k=1}^{\infty}\tilde\lambda_{m_k} \log\tilde\lambda_{m_k},$
which is finite since $\tilde{\lambda}_{m_{k}}=(1-2\epsilon)(1-r)r^{k-1}$ for some $r\in (2/3,1)$.
\ep
\medbreak

\acks We gratefully acknowledge A. Galves, R. Fern\'andez, K. Marton and S. Friedli for many discussion during the elaboration of this article. We thank the referees for the careful reading of the manuscript and for the comments, which improved the presentation of the article.

\bibliographystyle{jtbnew}
\bibliography{sandro_bibli.bib}

\begin{thebibliography}{30}
\expandafter\ifx\csname natexlab\endcsname\relax\def\natexlab#1{#1}\fi
\expandafter\ifx\csname url\endcsname\relax
  \def\url#1{\texttt{#1}}\fi
\expandafter\ifx\csname urlprefix\endcsname\relax\def\urlprefix{URL }\fi
\providecommand{\selectlanguage}[1]{\relax}

\bibitem[{Berger \emph{et~al.}(2005)Berger, Hoffman \&
  Sidoravicius}]{berger/hoffman/sidoravivius/2005}
\textsc{Berger, N., Hoffman, C. \& Sidoravicius, V.} (2005).
\newblock Nonuniqueness for specifications in $\ell^{2+\epsilon}$.
\newblock \emph{arXiv:math/0312344} .

\bibitem[{Bramson \& Kalikow(1993)}]{bramson/kalikow/1993}
\textsc{Bramson, M. \& Kalikow, S.} (1993).
\newblock Nonuniqueness in {$g$}-functions.
\newblock \emph{Israel J. Math.} \textbf{84}(1-2), 153--160.

\bibitem[{C\'enac \emph{et~al.}(2012)C\'enac, Chauvin, Paccaut \&
  Pouyanne}]{cenac/chauvin/paccaut/pouyanne/2010}
\textsc{C\'enac, P., Chauvin, B., Paccaut, F. \& Pouyanne, N.} (2012).
\newblock Variable length {M}arkov chains and dynamical sources.
\newblock \emph{S\'eminaire de Probabilit\'es XLIV, Lecture Notes in Math.}
  \textbf{2046}, 1--39.

\bibitem[{Chazottes \emph{et~al.}(2007)Chazottes, Collet, K{\"u}lske \&
  Redig}]{chazottes/collet/kulske/regig/2000}
\textsc{Chazottes, J.-R., Collet, P., K{\"u}lske, C. \& Redig, F.} (2007).
\newblock Concentration inequalities for random fields via coupling.
\newblock \emph{Probab. Theory Related Fields} \textbf{137}(1-2), 201--225.

\bibitem[{Comets \emph{et~al.}(2002)Comets, Fern{\'a}ndez \&
  Ferrari}]{comets/fernandez/ferrari/2002}
\textsc{Comets, F., Fern{\'a}ndez, R. \& Ferrari, P.~A.} (2002).
\newblock Processes with long memory: regenerative construction and perfect
  simulation.
\newblock \emph{Ann. Appl. Probab.} \textbf{12}(3), 921--943.

\bibitem[{De~Santis \& Piccioni(2012)}]{desantis/piccioni/2010}
\textsc{De~Santis, E. \& Piccioni, M.} (2012).
\newblock Backward coalescence times for perfect simulation of chains with
  infinite memory.
\newblock \emph{J. Appl. Probab.} \textbf{49}(2), 319--337.

\bibitem[{Doeblin \& Fortet(1937)}]{doeblin/fortet/1937}
\textsc{Doeblin, W. \& Fortet, R.} (1937).
\newblock Sur des cha\^{\i}nes \`a liaisons compl\`etes.
\newblock \emph{Bull. Soc. Math. France} \textbf{65}, 132--148.

\bibitem[{Fern{\'a}ndez \emph{et~al.}(2011)Fern{\'a}ndez, Gallo \&
  Maillard}]{fernandez/gallo/maillard/2011}
\textsc{Fern{\'a}ndez, R., Gallo, S. \& Maillard, G.} (2011).
\newblock Regular {$g$}-measures are not always {G}ibbsian.
\newblock \emph{Electron. Commun. Probab.} \textbf{16}, 732--740.

\bibitem[{Fern{\'a}ndez \& Maillard(2005)}]{fernandez/maillard/2005}
\textsc{Fern{\'a}ndez, R. \& Maillard, G.} (2005).
\newblock Chains with complete connections: general theory, uniqueness, loss of
  memory and mixing properties.
\newblock \emph{J. Stat. Phys.} \textbf{118}(3-4), 555--588.

\bibitem[{Friedli(2010)}]{friedli/2010}
\textsc{Friedli, S.} (2010).
\newblock A note on the {B}ramson-{K}alikow process.
\newblock \emph{Preprint
  \url{http://www.mat.ufmg.br/~sacha/textos/BK/plateaux.pdf}} .

\bibitem[{Gallo(2011)}]{gallo/2011}
\textsc{Gallo, S.} (2011).
\newblock Chains with unbounded variable length memory: perfect simulation and
  a visible regeneration scheme.
\newblock \emph{Adv. in Appl. Probab.} \textbf{43}(3), 735--759.

\bibitem[{Gallo \& Garcia(2011)}]{gallo/garcia/2011}
\textsc{Gallo, S. \& Garcia, N.~L.} (2011).
\newblock General context-tree-based approach to perfect simulation for chains
  of infinite order.
\newblock \emph{Submitted, arXiv:1103.2058v2} .

\bibitem[{Harris(1955)}]{harris/1955}
\textsc{Harris, T.~E.} (1955).
\newblock On chains of infinite order.
\newblock \emph{Pacific J. Math.} \textbf{5}, 707--724.

\bibitem[{Hulse(1991)}]{hulse/1991}
\textsc{Hulse, P.} (1991).
\newblock Uniqueness and ergodic properties of attractive {$g$}-measures.
\newblock \emph{Ergodic Theory Dynam. Systems} \textbf{11}(1), 65--77.

\bibitem[{Hulse(2006)}]{hulse/2006}
\textsc{Hulse, P.} (2006).
\newblock An example of non-unique {$g$}-measures.
\newblock \emph{Ergodic Theory Dynam. Systems} \textbf{26}(2), 439--445.

\bibitem[{Johansson \& {\"O}berg(2003)}]{johansson/oberg/2003}
\textsc{Johansson, A. \& {\"O}berg, A.} (2003).
\newblock Square summability of variations of {$g$}-functions and uniqueness of
  {$g$}-measures.
\newblock \emph{Math. Res. Lett.} \textbf{10}(5-6), 587--601.

\bibitem[{Kalikow(1990)}]{kalikow/1990}
\textsc{Kalikow, S.} (1990).
\newblock Random {M}arkov processes and uniform martingales.
\newblock \emph{Israel J. Math.} \textbf{71}(1), 33--54.

\bibitem[{Karlin(1953)}]{karlin/1953}
\textsc{Karlin, S.} (1953).
\newblock Some random walks arising in learning models. {I}.
\newblock \emph{Pacific J. Math.} \textbf{3}, 725--756.

\bibitem[{Keane(1972)}]{keane/1972}
\textsc{Keane, M.} (1972).
\newblock Strongly mixing {$g$}-measures.
\newblock \emph{Invent. Math.} \textbf{16}, 309--324.

\bibitem[{Lacroix(2000)}]{lacroix/2000}
\textsc{Lacroix, Y.} (2000).
\newblock A note on weak-{$\star$} perturbations of {$g$}-measures.
\newblock \emph{Sankhy\=a Ser. A} \textbf{62}(3), 331--338.

\bibitem[{Marton \& Shields(1994)}]{marton/shields/1994}
\textsc{Marton, K. \& Shields, P.~C.} (1994).
\newblock The positive-divergence and blowing-up properties.
\newblock \emph{Israel J. Math.} \textbf{86}(1-3), 331--348.

\bibitem[{McCullagh \& Nelder(1983)}]{mccullagh/nelder/1983}
\textsc{McCullagh, P. \& Nelder, J.~A.} (1983).
\newblock \emph{Generalized linear models}.
\newblock Monographs on Statistics and Applied Probability. London: Chapman \&
  Hall.

\bibitem[{Onicescu \& Mihoc(1935)}]{onicescu/mihoc/1935}
\textsc{Onicescu, O. \& Mihoc, G.} (1935).
\newblock Sur les cha\^ines de variables statistiques.
\newblock \emph{Bull. Sci. Math} \textbf{59}(2), 174--192.

\bibitem[{Preston(1976)}]{preston/1976}
\textsc{Preston, C.} (1976).
\newblock \emph{Random fields}.
\newblock Lecture Notes in Mathematics, Vol. 534. Berlin: Springer-Verlag.

\bibitem[{Propp \& Wilson(1996)}]{propp/wilson/1996}
\textsc{Propp, J.~G. \& Wilson, D.~B.} (1996).
\newblock Exact sampling with coupled {M}arkov chains and applications to
  statistical mechanics.
\newblock In: \emph{Proceedings of the {S}eventh {I}nternational {C}onference
  on {R}andom {S}tructures and {A}lgorithms ({A}tlanta, {GA}, 1995)}, vol.~9.

\bibitem[{Quas(1996)}]{quas/1996}
\textsc{Quas, A.~N.} (1996).
\newblock Non-ergodicity for {$C^1$} expanding maps and {$g$}-measures.
\newblock \emph{Ergodic Theory Dynam. Systems} \textbf{16}(3), 531--543.

\bibitem[{Rudolph(1982)}]{rudolph/1982}
\textsc{Rudolph, D.~J.} (1982).
\newblock A mixing {M}arkov chain with exponentially decaying return times is
  finitarily {B}ernoulli.
\newblock \emph{Ergodic Theory Dynam. Systems} \textbf{2}(1), 85--97.

\bibitem[{Shields(1996)}]{shields/1996}
\textsc{Shields, P.~C.} (1996).
\newblock \emph{The ergodic theory of discrete sample paths}, vol.~13 of
  \emph{Graduate Studies in Mathematics}.
\newblock Providence, RI: American Mathematical Society.

\bibitem[{Steif \& van~den Berg(1999)}]{steif/vanderberg/1999}
\textsc{Steif, J. \& van~den Berg, J.} (1999).
\newblock On the existence and nonexistence of finitary codings for a class of
  random fields.
\newblock \emph{Annals of Probability} \textbf{11}, 1501--1522.

\bibitem[{Stenflo(2001)}]{stenflo/2001}
\textsc{Stenflo, {\"O}.} (2001).
\newblock A note on a theorem of {K}arlin.
\newblock \emph{Statist. Probab. Lett.} \textbf{54}(2), 183--187.

\end{thebibliography}

\end{document}